\documentclass[10pt]{amsart}
\pdfoutput=1
\usepackage{url}
\usepackage{graphicx}
\usepackage{caption}

\usepackage[margin=1in]{geometry}

\DeclareGraphicsExtensions{.png,.eps}

\newtheorem{theorem}{Theorem}

\theoremstyle{definition}

\newcommand{\incentive}{\varphi}

\makeatletter
\g@addto@macro{\endabstract}{\@setabstract}
\newcommand{\authorfootnotes}{\renewcommand\thefootnote{\@fnsymbol\c@footnote}}%
\makeatother
\begin{document}


\begin{center}
  \LARGE Entropic Equilibria Selection of Stationary Extrema in Finite Populations
   \par \bigskip

  \normalsize
  \authorfootnotes
  Marc Harper\footnote{corresponding author, email: marc.harper@gmail.com}, Dashiell Fryer\textsuperscript{1}\par \bigskip

  \textsuperscript{1}Department of Mathematics, Pomona College\par \bigskip
  
\end{center}

\begin{abstract}
We propose the entropy of random Markov trajectories originating and terminating at a state as a measure of the stability of a state of a Markov process. These entropies can be computed in terms of the entropy rates and stationary distributions of Markov processes. We apply this definition of stability to local maxima and minima of the stationary distribution of the Moran process with mutation and show that variations in population size, mutation rate, and strength of selection all affect the stability of the stationary extrema.
\end{abstract}

\section{Introduction}

This work is motivated by the stationary stability theorem \cite{harper2013stationary}, which characterizes local maxima and minima of the stationary distribution of the Moran process with mutation in terms of evolutionary stability. Specifically, the theorem says that for sufficiently large populations, the local maxima and minima of the stationary distribution satisfy a selective-mutative equilibria criterion that generalizes the celebrated notion of evolutionary stability \cite{smith1982evolution}. This means that the stationary distribution encodes the usual information about evolutionary stability. Precisely which equilibria are favored (i.e. are maxima or minima) is a natural question and depends on the choice of various parameters, such as the mutation rate $\mu$, the strength of selection $\beta$, and the population size $N$.

We propose the random trajectory entropy (RTE) of paths originating and terminating at a state as a measure of stability of the state \cite{ekroot1993entropy} \cite{kafsi2013entropy}. This is an information-theoretic quantity that is easily computable from the entropy rate and stationary distribution of a process, and varies continuously with the critical evolutionary parameters (as does the stationary distribution). We will see that RTE captures the behavior of the Moran process with mutation intuitively, leading to a simple method for equilibrium selection for finite populations, generally a significant problem in evolutionary game theory \cite{samuelson1998evolutionary} \cite{harsanyi1988general}.

\section{Stationary Distributions, Entropy Rates, and Random Trajectory Entropies}

Our first goal is to establish the random trajectory entropy (RTE) of a state as a measure of instability of the state. We will be particularly concerned with the local and global extrema of the stationary distribution, shown in \cite{harper2013stationary} to have a close connection with evolutionary stability.  

The stationary distribution of a Markov process gives the probability that the process will be in each state in the long run \cite{hordijk1988insensitive}. As such it is a fundamental convergence concept for Markov processes. We take the weighted graph viewpoint of Markov processes on a finite set of states $V$. Let the transition probabilities be given by a function $T:V \times V \to [0,1]$ (viewed as a matrix or a function), and the stationary distribution by a function $s: V \to [0,1]$ (appropriately normalized to a probability distribution).  We assume throughout that all processes are irreducible (there is a path between any two states) and have unique stationary distributions.Let $V' \subset V$ and define a stationary maximum of $V'$ to be a state $v \in V'$ such that $s(v') < s(v)$ for all $v' \in V' \setminus {v}$. Then we have a local maximum $v$ if the set $V'$ is the set of neighboring states of $v$ and a global maximum if $V' = V$; similarly for minima.


Although the stationary distribution of a process is often quite useful, it does not tell the full story of the process. While the stationary distribution gives the long run occupancy of any particular state, it does not explain how much the process moves among states, and so gives an incomplete description of the dynamic stability of a state. An entropy rate is a generalization of Shannon entropy to Markov processes and are commonly described as the \emph{inherent randomness} or \emph{information content} of a process \cite{ekroot1993entropy}. The entropy rate of a process encodes both long term and short term information about the process, defined for a process $X$ as follows:
\begin{equation} \label{entropy_rate}
 H(X) = - \sum_{i, j}{s(v_i) T(v_i, v_j) \log T(v_i, v_j)}
\end{equation}
The entropy rate is a value attached a process rather individual states. We need a quantity associated to both the process and the individual states that can discriminate between equilibria.

Following \cite{ekroot1993entropy}, define the probability of a trajectory $V: v_0 \to v_1 \to \cdots \to v_k$ with no intermittent state being $v_k$ as the product of the transitions along the path
\[ Pr(V) = T(v_0, v_1) T(v_1, v_2) \cdots T(v_{k-1}, v_k).\]
Since the process is irreducible, we have that the sum over all possible such trajectories from $v_0$ to $v_k$ is one, forming a probability distribution. Let $\mathcal{T}(v_0, v_k)$ be the set of all such paths and define the \emph{random trajectory entropy} (RTE) from $v_0$ to $v_k$ to be the entropy of the probability distribution on $\mathcal{T}(v_0, v_k)$, i.e.
\[ H_{v_0 v_k} = - \sum_{v \in \mathcal{T}(v_0, v_k)}{Pr(v) \log Pr(v)}.\]
It was shown in \cite{ekroot1993entropy} (Theorem 1, p.1419) that when the starting and ending states are the same, the entropy of the random trajectory is determined by the entropy rate and the stationary probability
\[ H_{v} := H_{v v} = \frac{H(X)}{s(v)}. \]
From this we have immediately have the following theorem characterizing local and global extrema of the stationary distribution.
\begin{theorem}\label{theorem_entropy_rate}
For an irreducible Markov process with stationary distribution $s$, a state $s$ is a local (resp. global) maximum (resp. minimum) if and only if the RTE $H_{v}$ is a local (resp. global) minimum (resp. maximum).
\end{theorem}

Furthermore we now recognize the random trajectory entropy as a measure of stationary instability of a state, which we can now use to compare and select equilibria for the same process and for closely related processes. Intuitively, a smaller RTE means that trajectories tend to stay near a local maxima, i.e that random walks tend to be short, which is a way of saying that the state is \emph{stable}. (Note that $1/s(v)$ is the expected number of steps it takes to return to $v$.)

\section{Applications}

We now consider several explicit examples of finite population processes.

\subsection{Moran Process with Mutation}

For the Moran process with mutation we use a special case of the formulation \cite{harper2013stationary}; see also \cite{fudenberg2004stochastic}, \cite{claussen2005non}, and \cite{moran1962statistical}. Let a population be composed of $n$ types $A_1, \ldots A_n$ of size $N$ with $a_i$ individuals of type $A_i$ so that $N = a_1 + \cdots + a_n$. We will denote a population state by the tuple $a = (a_1, \ldots, a_n)$ and the population distribution by $\bar{a} = a / N$. We assume the existence of a fitness landscape $f$ where $f_i(\bar{a})$ gives the fitness of type $A_i$; typically $f(\bar{a}) = G \bar{a}$ for some game matrix $G$. (See \cite{hofbauer2003evolutionary} \cite{weibull1997evolutionary} and \cite{hofbauer1998evolutionary} for general references on evolutionary games). Define a matrix of mutations $M$ where $0 \leq M_{i j} \leq 1$ may be a function of the population state for our most general results, but we will typically assume in examples that for some constant value $\mu$, the mutation matrix takes the form $M_{i j} = \mu / (n-1)$ for $i \neq j$ and $M_{i i} = 1 - \mu$. A typical mutation rate is $\mu \approx 1/N$.

The Moran process with mutation is a Markov process on the population states defined by the following transition probabilities, corresponding to a birth-death process where birth is fitness-proportionate with mutation and death is uniformly random. To define the adjacent population states, let $i_{\alpha \beta}$ be the vector that is 1 at index $\alpha$, -1 at index $\beta$, and zero otherwise, with the convention that $i_{\alpha \alpha}$ is the zero vector of length $n$. Every adjacent state of state $a$ for the Moran process is of the form $a + i_{\alpha \beta}$ for some $1 \leq \alpha, \beta \leq n$. At a population state $a$ we choose an individual of type $A_i$ to reproduce proportionally to its fitness, allowing for mutation of the new individual as given by the mutation probabilities. The distribution of fitness proportionate selection probabilities is given by $p(\bar{a}) = M(\bar{a}) \bar{\incentive}(\bar{a})$; explicitly, the $i$-th component is
\begin{equation}
p_i(\bar{a}) = \frac{\sum_{k=1}^{n}{\incentive_k(\bar{a}) M_{k i} }}{\sum_{k=1}^{n}{\incentive_k(\bar{a})}}
\label{fitness-proportionate-reproduction} 
\end{equation}
where the function $\incentive(\bar{a}) = \bar{a}_i f_i(\bar{a})$. We also randomly choose an individual to be replaced, just as in the Moran process. This yields the transition probabilities
\begin{align}
T_{a}^{a + i_{\alpha, \beta}} &= p_{\alpha}(\bar{a}) \bar{a}_{\beta} \qquad \text{for $\alpha \neq \beta$} \notag \\
T_{a}^{a} &= 1 - \sum_{b \text{ adj } a, b\neq a}{T_{a}^{b}}
\label{moran_process}
\end{align}

We will also utilize a variant incorporating a strength of selection term $\beta$ called Fermi selection \cite{traulsen2009stochastic}:
\[ \incentive(\bar{a}) = \bar{a}_i e^{\beta f_i(\bar{a})}.\] 
For our examples we will restrict our attention to processes defined by $X = X(N, n, \mu, \incentive)$. Several explicit examples of stationary distributions for Moran processes with mutation are given in \cite{claussen2005non} \cite{harper2014inherent}. The entropy rate of the Moran process with mutation was computed in \cite{harper2014inherent} ($n=2$) and \cite{harper2014entropy} ($n>2$) along with the development of a number of theoretical results. For our purposes the explicit values of the entropy rate are not needed. Generally as $\mu \to 0$, the entropy rate also goes to zero, and attains its maximum as $N \to \infty$ for the neutral fitness landscape (with e.g. mutations $\mu = 1/N$). The entropy rate of is bounded below by zero and above by $\frac{2n-1}{n} \log{n}$ \cite{harper2014entropy}. The RTE is bounded below by the entropy rate, justifying the description of the entropy rate as the inherent randomness of a process.

\subsection{Comparison of Equilibria of a Single Process}

Since the entropy rate is associated to the entire process, for two different states $i$ and $j$ we have that $H_{i} = H(X)/s(i)$ and $H_{j} = H(X) / s(j)$, so if $H(X) \neq 0$ we need only consider the values of the stationary process in this case to compare the equilibria as $H_{j} / H_{i} = s(i) / s(j)$.

\subsubsection{Small mutation limit} Consider the special case in which the rate of mutation parameter $\mu \to 0$ in a population of two types. Then we have $\lim_{\mu \to 0}{H(X)} = 0$ \cite{harper2014entropy}. In this case the stationary distribution becomes a delta distribution on the corner states. For population of two types $A$ and $B$, we can express the limiting stationary distribution in terms of the fixation probabilities of the two types $\rho_A$ and $\rho_B$ \cite{fudenberg2004stochastic}:
\[ \lim_{\mu \to 0} s(0, N) = \frac{\rho_B}{\rho_A + \rho_B} \quad \text{and} \quad \lim_{\mu \to 0} s(N, 0) = \frac{\rho_A}{\rho_A + \rho_B}\]
Hence we have that
\[ \lim_{\mu \to 0}\frac{H_{(0, N)}}{H_{(N, 0)}} = \lim_{\mu \to 0} \frac{s(N, 0)}{s(0, N)} = \frac{\rho_B}{\rho_A}\]
In other words, the state with the type having greater fixation probability is more stable. For the classical Moran process with game matrix
$G = \left( \begin{smallmatrix}
 r & r\\
 1 & 1\\
\end{smallmatrix} \right) $
we have that (assuming $r \neq 1$): $\rho_A = (1 - r^{-1}) / (1 - r^{-N})$ and $\rho_B = (1 - r^{1-N}) / (1 - r^{-N})$, which gives
\[ \lim_{\mu \to 0} \frac{H_{(0, N)}}{H_{(N, 0)}} = \frac{\rho_B}{\rho_A} = \frac{1 - r^{1-N}}{1 - r^{-1}}\]
As expected, whether $r>1$ determines which equilibrium is favored. If $r = 1$, $\rho_A = 1/N = \rho_B$ and the RTEs are equal.

\subsubsection{Large Populations and Neutral Landscapes}

For arbitrarily many types, the stationary distribution for the neutral fitness landscape (matrix of all ones) and any mutation rate $\mu$ can be analytically computed. For large $N$, the neutral landscape attains the maximum entropy rate, so for large populations a sufficient condition for a state for a non-neutral landscape to be more stable than the same state for the neutral landscape is simply to have a larger stationary probability \cite{harper2014entropy}. For non-neutral landscapes, the large population limit need not maximize the entropy rate \cite{harper2014entropy}.

\subsection{Comparison of Equilibria for Separate Processes on the Same States}

\begin{figure}[H]
 \includegraphics[width=\textwidth]{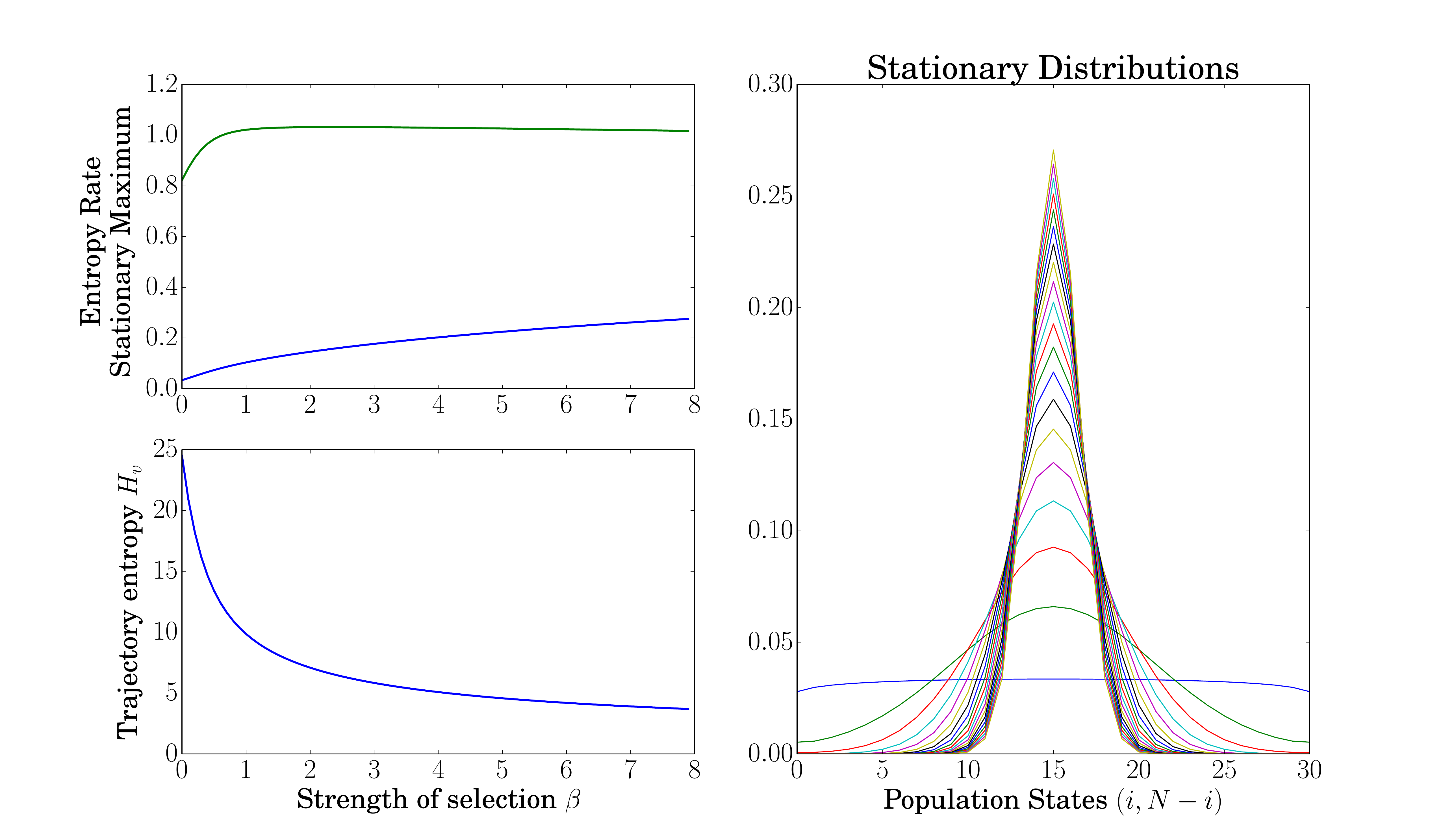}
 \caption{Right: Stationary distributions for Hawk-Dove landscapes for varying strength of selection $\beta \in [0, 10]$, $N=30$, $\mu = 1/N$.  Upper Left: As $\beta$ increases, so does the stationary probability (blue, lower curve) of the maxima at $(15, 15)$. The entropy rate (green, upper) is not monotonically increasing in $\beta$. Lower Left: Nevertheless, as $\beta$ increases, the random trajectory entropy decreases monotonically as expected intuitively. More intense selection yields greater stability at the maximum.}
 \label{figure_RTE_HD}
\end{figure}

Two instances of the Moran process with mutation can have the same stationary maximum state but different entropy rates. Consider the one-parameter family corresponding to a Hawk-Dove matrix 
\begin{equation} \label{hd_matrix}
 G = \left( \begin{smallmatrix}
 1 & 2\\
 2 & 1\\
\end{smallmatrix} \right),
\end{equation}
and transition probabilities defined by Fermi selection. For convenience fix a population size $N \geq 10$ and $N$ even, and let the rate of mutation be $\mu = 1/N$. Then we have that $(N/2, N/2)$ is the unique stationary maximum \cite{harper2013stationary}. As $\beta$ increases, the stationary probability at the maximum increases more quickly than the entropy rate (which is not monotonic in this case). The net result is that the random trajectory entropy is decreasing as a function of $\beta$ and that the stationary maximum is ``more stable'' (as would is expected for greater strengths of selection). See Figure \ref{figure_RTE_HD}. For other two player games the situation is analogous, e.g. for the Coordination game, the interior equilibrium RTE is decreasing as a function of $\beta$. For both games, we have the intuitive result that the stability (measured by the RTE) of the extrema varies monotonically with the strength of selection.

\begin{figure}
 \includegraphics[width=\textwidth]{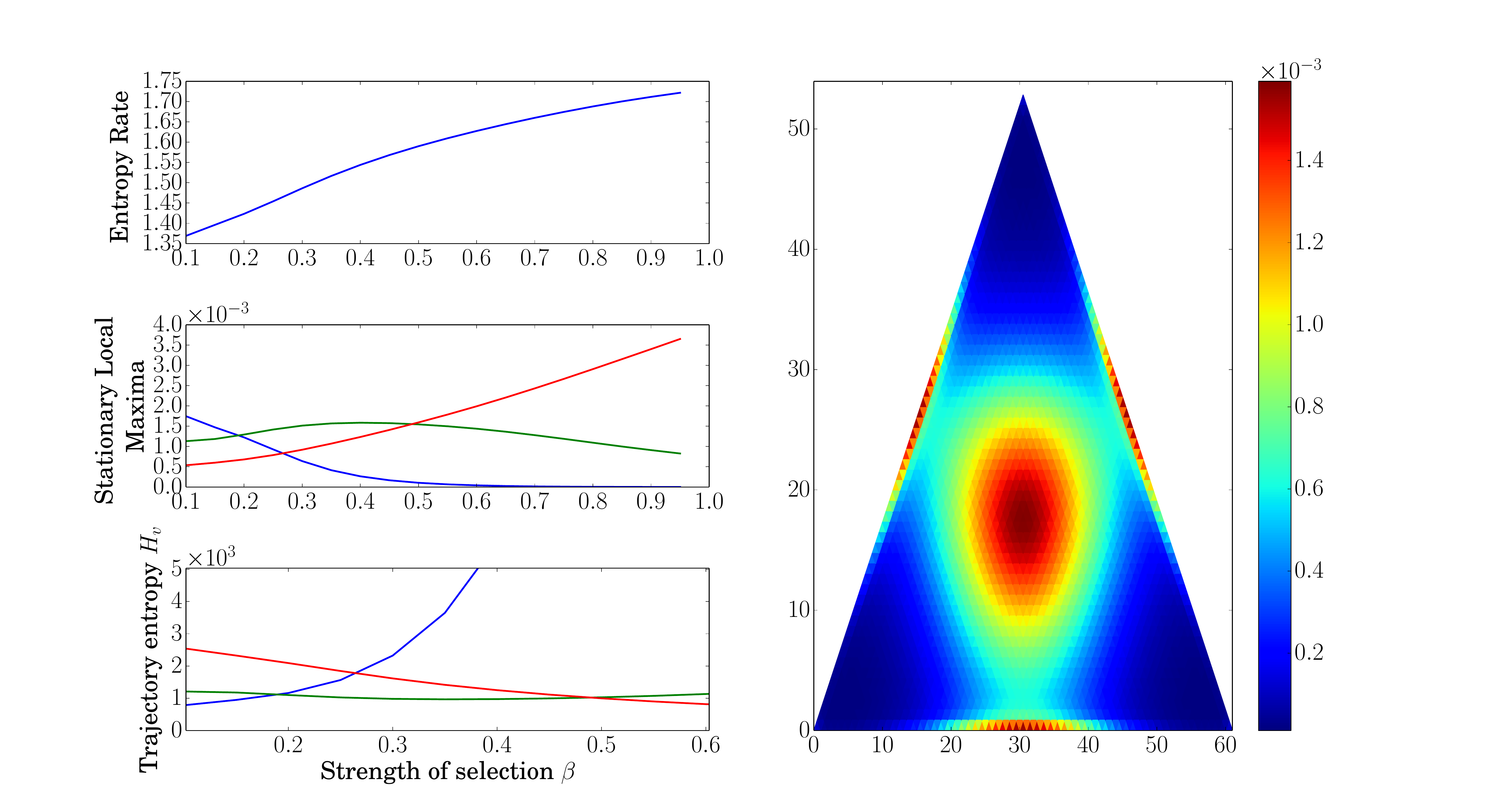}
 \caption{ This $n=3$ player example for landscape defined by the matrix \ref{3x3matrix}, $N=60$, $\mu=1/N$ has multiple local stationary extrema, at the center of the simplex, on the centers of the boundary simplices, and on the corners of the simplex. Left Top: The entropy rate of the process is given as a function of the strength of selection $\beta$. Left Center: As $\beta$ increases, the stationary probability of each extrema changes. The curves correspond as follows: Blue $(N,0,0)$, Green $(N/2, N/2, 0)$, Red $(N/3, N/3, N/3)$. (Symmetric permutations of these states are also extrema and have the same probabilities.) As the strength of selection increases, more stationary probability is concentrated on the central extrema. Lower Left: As $\beta$ increases, the trajectory entropy of the boundary extrema increases while decreasing for the central extrema, showing that strength of selection affects the stability of the equilibria. Which of the equilibria is most stable depends on the value of $\beta$. Right: Stationary distribution for $\beta = 0.5$.}
 \label{random_entropy_3}
\end{figure}

\begin{figure}
 \includegraphics[width=0.8\textwidth]{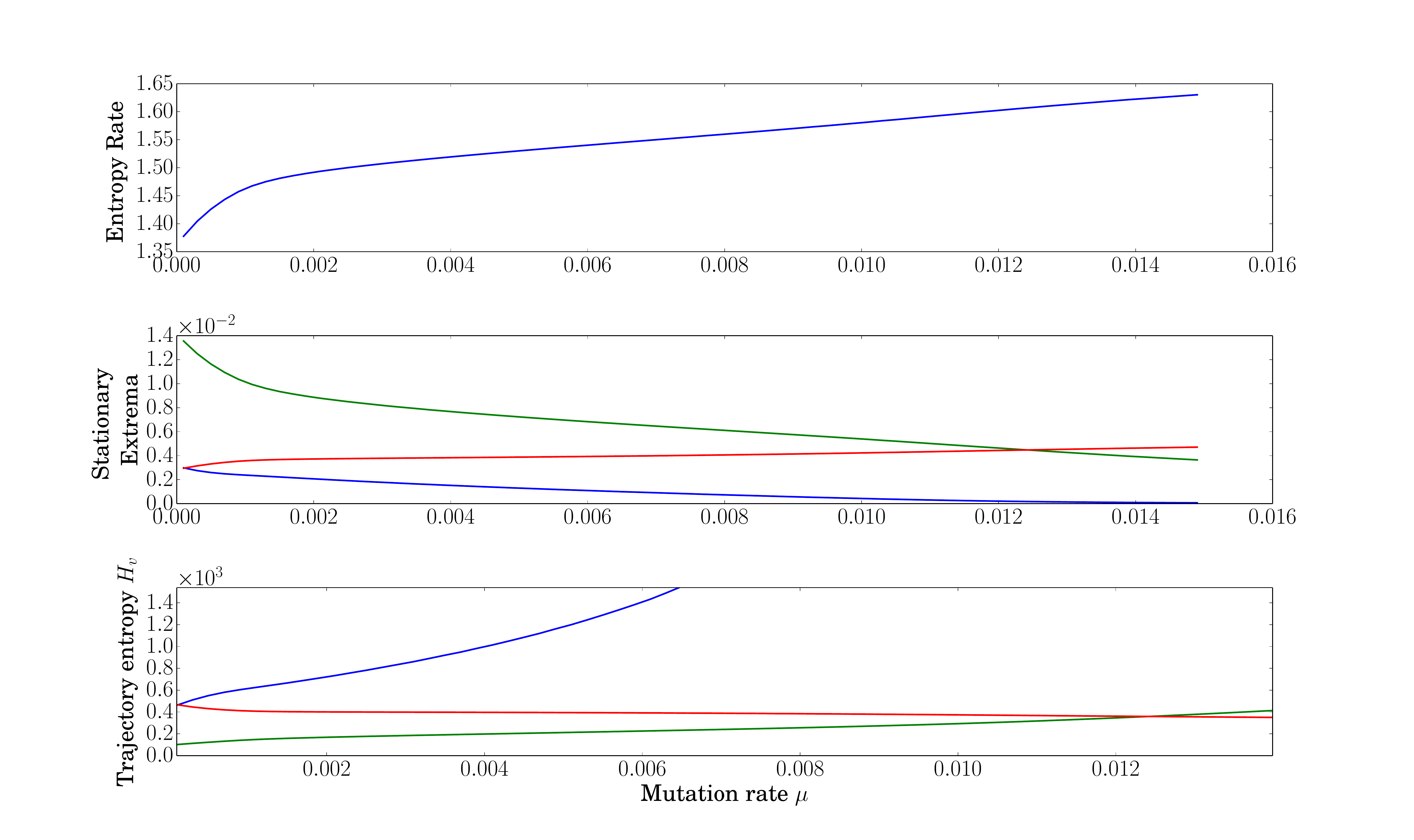}
 \caption{Entropy Rate (Upper Left), Stationary probabilities (Center Left), and Trajectory entropies (Lower Left) for a process with $N=42$, $\beta=1$, landscape defined by Matrix (\ref{3x3matrix}) and varying rate of mutation $\mu$. Just as for the strength of selection $\beta$ in Figure \ref{random_entropy_3}, the value of $\mu$ can determine which of the equilibria is most stable. As $\mu \to 0$ the corner states are favored. As $\mu$ increases, the interior equilibrium becomes more stable.}
 \label{figure_trajectory_mu}
\end{figure}

We now consider multiple examples for the landscape derived from the three-type game matrix:
\begin{equation} \label{3x3matrix}
 G = \left( \begin{smallmatrix}
 0 & 1 & 1\\
 1 & 0 & 1\\
 1 & 1 & 0\\
\end{smallmatrix} \right) 
\end{equation}
This landscape typically has several local extrema. Let the population size $N'=6N$. Then we have extrema at the simplex corners $(6N,0,0), (0,6N,0), (0,0,6N)$, center of the boundary simplices $(3N, 3N, 0), (3N, 0, 3N), (0, 3N, 3N)$, and the center $(2N, 2N, 2N)$. Varying either $\mu$ (Figure \ref{figure_trajectory_mu}) or $\beta$ (Figure \ref{random_entropy_3}) changes which equilibria have the smallest RTE. As $\mu$ increases, stationary probability moves from the corner points of the simplex to the midpoints of the boundary simplices and also toward the center. Similarly for the strength of selection $\beta$. In both cases the rate of change of the stationary extrema dominates since the entropy rate varies slowly.

\begin{figure}
 \includegraphics[width=0.5\textwidth]{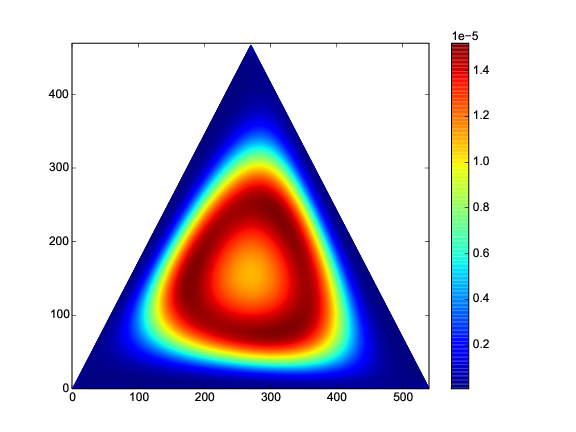}
 \caption{Stationary distribution for a RSP landscape (matrix below) and large population $N=540$, $\beta=1$, $\mu = 3/(2N)$. There are apparent cycles of constant stationary probability and hence constant RTE. This is analogous to the concentric cycles of the replicator equation \cite{hofbauer1998evolutionary}. The central state is a local stationary maximum and has a large trajectory entropy.}
 \label{figure_rsp}
\end{figure}

Though we have focused on equilibria, the stationary distributions of finite population games can exhibit a variety of complex dynamical behaviors such as depicted in Figure \ref{figure_rsp}. Consider the rock-paper-scissors landscape given by the matrix
\[G = \left( \begin{smallmatrix}
 0 & 1 & -1\\
 -1 & 0 & 1\\
 1 & -1 & 0\\
 \end{smallmatrix} \right).\]
 For some parameter choices, RPS landscapes produce an interesting stationary distribution with discretized cycles of constant trajectory entropy, analogous to the concentric cycles for the replicator equation and the fact that the relative entropy is a constant of motion of the replicator equation \cite{weibull1997evolutionary}. Assuming symmetry of the cycle (a large population of size divisible by 3 seems to suffice to yield approximate symmetry), no value on any cycle is a local maximum and the values on the maximal cycle are all global maxima. In this case the stationary stability theorem does not apply to the cycles (but still applies to the local minimum in the center of the simplex).

\subsection{Comparison of Equilibria for Process with Varying Population Size}

\begin{figure}
 \includegraphics[width=0.8\textwidth]{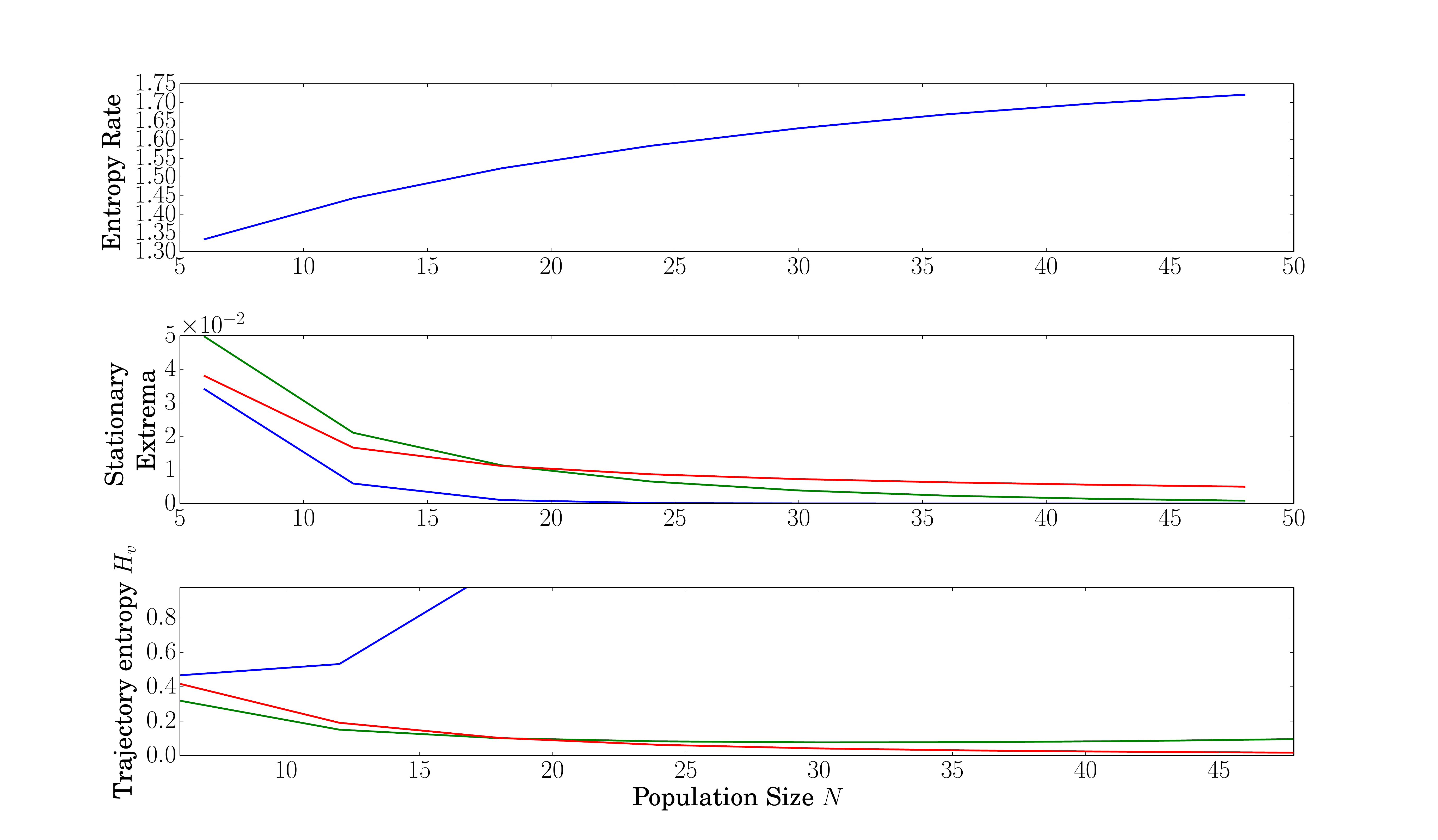}
 \caption{Entropy Rate (Upper), Stationary probabilities (Center), and Random Trajectory entropies (Lower) for a process with $\beta=1$, landscape defined by \ref{3x3matrix}, varying population size $N$, and $\mu=1/N$. Just as for the strength of selection $\beta$ in Figure \ref{random_entropy_3} and $\mu$ in Figure \ref{figure_trajectory_mu}, the population size $N$ can determine which of the equilibria is most stable. The trajectory entropies have been scaled by the number of states of the process, $\binom{N+3}{3}$.}
 \label{figure_trajectory_N}
\end{figure}

For the final example we consider the effect of altering the population size the population size $N$. In this case the underlying state spaces are different even though the equilibria are generally same) for large enough $N$. For the same number of types $n$ the entropy rate has the same upper bound (though the entropy rate typically increases with $N$), and so to enable a fair comparison we normalize by the number of states (since the stationary distribution is spread out over a variable number of states). In general, the number of states is $\binom{N+n-1}{n}$. As for both $\beta$ and $\mu$, varying the population size $N$ changes the favored equilibrium. See Figure \ref{figure_trajectory_N}. We note, however, that the RTEs are increasing in $N$, and so the issue of normalization is critical to the comparison of equilibria for processes with different population sizes.

\section{Discussion}

We have proposed random trajectory entropy as a measure of stability of states of finite Markov processes and considered several examples from finite population biology. Variations of fundamental evolutionary parameters alters the stability of equilibria and agrees with intuitive expectations. In particular, stability is closely tracked by stationary probability in several example population dynamics. We did not consider RTEs for paths that occur originate and terminate at different states but it is reasonable to expect that e.g. that a local stationary maxima will have smaller RTE in some neighborhood, and similarly for local minima.

All computations were performed with open source code available at \url{https://github.com/marcharper/stationary}. This package can compute exact stationary distributions and entropy rates for reversible processes and approximate solutions for all other cases mentioned in this manuscript. All plots created with \emph{matplotlib} \cite{Hunter:2007}.

\bibliographystyle{unsrt}
\bibliography{ref}

\begin{thebibliography}{10}

\bibitem{harper2013stationary}
Marc Harper and Dashiell Fryer.
\newblock Stationary stability for evolutionary dynamics in finite populations.
\newblock {\em arXiv preprint arXiv:1311.0941}, 2013.

\bibitem{smith1982evolution}
John~Maynard Smith.
\newblock {\em Evolution and the Theory of Games}.
\newblock Cambridge university press, 1982.

\bibitem{ekroot1993entropy}
Laura Ekroot and Thomas~M Cover.
\newblock The entropy of markov trajectories.
\newblock {\em Information Theory, IEEE Transactions on}, 39(4):1418--1421,
  1993.

\bibitem{kafsi2013entropy}
Mohamed Kafsi, Matthias Grossglauser, and Patrick Thiran.
\newblock The entropy of conditional markov trajectories.
\newblock {\em Information Theory, IEEE Transactions on}, 59(9):5577--5583,
  2013.

\bibitem{samuelson1998evolutionary}
Larry Samuelson.
\newblock {\em Evolutionary games and equilibrium selection}, volume~1.
\newblock Mit Press, 1998.

\bibitem{harsanyi1988general}
John~C Harsanyi and Reinhard Selten.
\newblock A general theory of equilibrium selection in games.
\newblock {\em MIT Press Books}, 1, 1988.

\bibitem{hordijk1988insensitive}
Arie Hordijk and Ad~Ridder.
\newblock Insensitive bounds for the stationary distribution of non-reversible
  markov chains.
\newblock {\em Journal of applied probability}, pages 9--20, 1988.

\bibitem{fudenberg2004stochastic}
Drew Fudenberg, Lorens Imhof, Martin~A Nowak, and Christine Taylor.
\newblock Stochastic evolution as a generalized moran process.
\newblock {\em Unpublished manuscript}, 2004.

\bibitem{claussen2005non}
Jens~Christian Claussen and Arne Traulsen.
\newblock Non-gaussian fluctuations arising from finite populations: Exact
  results for the evolutionary moran process.
\newblock {\em Physical Review E}, 71(2):025101, 2005.

\bibitem{moran1962statistical}
Patrick Alfred~Pierce Moran et~al.
\newblock The statistical processes of evolutionary theory.
\newblock {\em The statistical processes of evolutionary theory.}, 1962.

\bibitem{hofbauer2003evolutionary}
Josef Hofbauer and Karl Sigmund.
\newblock Evolutionary game dynamics.
\newblock {\em Bulletin of the American Mathematical Society}, 40(4):479--519,
  2003.

\bibitem{weibull1997evolutionary}
J{\"o}rgen~W Weibull.
\newblock {\em Evolutionary game theory}.
\newblock MIT press, 1997.

\bibitem{hofbauer1998evolutionary}
Josef Hofbauer and Karl Sigmund.
\newblock {\em Evolutionary games and population dynamics}.
\newblock Cambridge University Press, 1998.

\bibitem{traulsen2009stochastic}
Arne Traulsen and Christoph Hauert.
\newblock Stochastic evolutionary game dynamics.
\newblock {\em Reviews of nonlinear dynamics and complexity}, 2:25--61, 2009.

\bibitem{harper2014inherent}
Marc Harper.
\newblock Inherent randomness of evolving populations.
\newblock {\em Physical Review E}, 89(3):032709, 2014.

\bibitem{harper2014entropy}
Marc Harper.
\newblock Entropy rates of the multidimensional moran processes and
  generalizations.
\newblock {\em arXiv preprint arXiv:1401.2713}, 2014.

\bibitem{Hunter:2007}
J.~D. Hunter.
\newblock Matplotlib: A 2d graphics environment.
\newblock {\em Computing In Science \& Engineering}, 9(3):90--95, 2007.

\end{thebibliography}

\end{document}